\documentclass[12pt]{amsart}
\usepackage{amssymb,amsmath,latexsym,amscd,amsfonts,amsrefs,color}
\usepackage{graphics}
\usepackage{epsfig}
\usepackage{psfrag}
\usepackage{enumerate}
\def\ignore #1 {}

\newtheorem{thm}{Theorem}
\newtheorem{prop}[thm]{Proposition}
\newtheorem{lem}[thm]{Lemma}
\newtheorem{cor}[thm]{Corollary}

\theoremstyle{definition}
\newtheorem{dfn}[thm]{Definition}
\newtheorem{rem}[thm]{Remark}
\newtheorem{ex}[thm]{Example}

\def\hpic #1 #2 {\mbox{$\begin{array}[c]{l} \epsfig{file=#1,height=#2} \end{arr\
ay}$}}
\def\vpic #1 #2 {\mbox{$\begin{array}[c]{l} \epsfig{file=#1,width=#2} \end{arra\
y}$}}

\begin{document}

\title{Asymptotic Expansion of Warlimont Functions on Wright Semigroups}

\author{Marco Aldi and Hanqiu Tan}

\begin{abstract}
We calculate full asymptotic expansions of prime-independent multiplicative functions on additive arithmetic semigroup that satisfy a strong form of Knopfmacher's axioms. When applied to the semigroup of unlabeled graphs, our method yields detailed asymptotic information on how graphs decompose into connected components. As a second class of examples, we discuss polynomials in several variables over a finite field.
\end{abstract}


\maketitle

\section{Introduction}

Let $G_n$ be the number of unlabeled graphs with $n$ vertices and let $G_n^+$ be the number of connected unlabeled graphs with $n$ vertices. Using the fact that the sequences $\{G_n\}$ and $\{G_n^+\}$ are related by the identity 
\begin{equation}\label{eq:00}
\sum_{n=0}^\infty G_n x^n = \prod_{m=1}^\infty (1-x^{m})^{-G_m^+}\,,
\end{equation}
Wright \cite{wright67} proved that $G_n\sim G_n^+$ i.e.\ almost all graphs are connected. Armed with a full asymptotic expansion for $G_n$ \cite{wright69}, Wright further improved this result by constructing \cite{wright70} a sequence $\{\omega_s\}$ of polynomials such that
\begin{equation}\label{eq:000}
G_n^+=G_n+\sum_{s=1}^{R-1} \omega_s(n)G_{n-s}+O(n^RG_{n-R})
\end{equation}
for all positive integers $R$.

In the context of abstract analytic number theory \cite{K2}, Knopfmacher \cite{K} observed that \eqref{eq:00} is a particular case of an Euler product type of identity that holds for arbitrary additive arithmetical semigroups and the methods of \cite{wright67} can be used to study the distribution of certain arithmetical functions on additive arithmetical semigroups in which almost all elements are prime. For instance, if $d_2$ is the {\it divisor function} that to each unlabeled graph $g$ assigns the number of ways to write $g$ as a disjoint union of an ordered pair of graphs then
\begin{equation}\label{eq:0000}
\lim_{n\to \infty} \frac{1}{G_n}\sum d_2(g)=2\qquad\textrm{ and }\qquad
\lim_{n\to \infty} \frac{1}{G_n}\sum (d_2(g)-2)^2=0\,,
\end{equation}
where both sums are taken over all graphs $g$ with $n$ vertices. 

The goal of the present paper is to investigate Knopfmacher's suggestion \cite{K} that restricting to  arithmetical semigroups in which the total number of elements is related to the number of prime elements by a formula analogous to \eqref{eq:000} might lead to a strengthening of \eqref{eq:0000}. To illustrate our results with an example, consider again the divisor function $d_2$ on the semigroup of graphs. We prove that for every positive integer $M$, there exists a sequence $\{\tau_s(n)\}$ of polynomials such that
\begin{equation}\label{eq:00000}
\frac{1}{G_n}\sum (d_2(g)-2)^M = 2^M\sum_{s=1}^{R-1} \tau_s(n) 2^{-sn} + O(n^{2R-1}2^{-Rn})
\end{equation}
for every integer $R$. Clearly, \eqref{eq:0000} can be recovered by setting $M=1$ and $M=2$ in \eqref{eq:00000} and taking the limit as $n\to \infty$. More generally, we show that \eqref{eq:00000} is a particular case of a formula that holds if $d_2$ is replaced by an arbitrary {\it Warlimont function} i.e.\ a multiplicative prime-independent function whose restriction to power of primes grows in a prescribed way. Even more generally, the semigroup of graphs can be replaced by any {\it Wright semigroup} which we define to be an additive arithmetical semigroup subject to a growth condition introduced in \cite{wright70}. Examples of Wright semigroups include the semigroup of unlabeled graphs with an even number of edges and the semigroup of polynomials in at least two variables over a finite field.

The paper is organized as follows. Section \ref{sec:2} contains the main technical results used in the rest of the paper. We work with triples of sequences related by a generalization of \eqref{eq:00} that were introduced in \cite{warlimont}. The main result is Theorem \ref{thm:5} which can be thought of as a generalization of \cite{wright70}, modeled after the way in which \cite{warlimont} generalizes \cite{wright67}. In Section \ref{sec:3}, after introducing the key notions of Wright semigroup and of Warlimont function, provide asymptotic formulas for moments of Warlimont functions in terms of the number of elements of given degree in the underlying (not necessarily Wright) semigroup. In the special case of Wright semigroup, we construct full asymptotic expansions generalizing \eqref{eq:00000}. We illustrate our results in Section 4 by calculating asymptotic expansion of some of the arithmetical functions considered in \cite{K} on three examples of Wright semigroups: the semigroup of all unlabeled graphs, the semigroup of unlabeled graphs with an even number of edges and the semigroup of non-zero polynomials (up to scaling) in at least two variables over a finite field.

\bigskip\noindent \textbf{Acknowledgments:} Part of this work was carried out during the Summer of 2016 at VCU and supported by a UROP Summer Research Fellowship.

\section{Warlimont Triples}\label{sec:2}

\begin{dfn}
A {\it Warlimont triple} is a triple $(\{T_n\},\{t_n\},\{a_n\})$ of sequences of non-negative real numbers related by the following identity of formal power series
\begin{equation}\label{eq:0}
\sum_{n=0}^\infty T_n x^n = \prod_{m=1}^\infty \left(\sum_{k=0}^\infty a_k x^{km}\right)^{t_m} 
\end{equation}
and such that 
\begin{enumerate}[i)]
\item $T_0=a_0=1$;
\item $a_1>0$;
\item $t_m\in \mathbb Z$ for all $m$ and $t_m>0$ for all but finitely many $m$.
\end{enumerate}
\end{dfn}

\begin{lem}\label{lem:2}
Let $(\{T_n\},\{t_n\},\{a_n\})$ be a Warlimont triple and consider the sequences $\{v_n\}$, $\{\beta_n\}$ and $\{b_n\}$ defined the recursion formulas
\begin{align}
v_n&=T_n-\sum_{s=1}^{n-1}\frac{s}{n} v_sT_{n-s}\label{eq:1}\\
\beta_n&=-\sum_{s=0}^{n-1}\beta_sT_{n-s}\label{eq:2}\\
b_n&=na_n - \sum_{s=1}^{n-1} b_sa_{n-s}\label{eq:3}  
\end{align}
with initial conditions $v_1=T_1$, $\beta_0=1$, $b_1=a_1$. Then for all $n$
\begin{enumerate}[1)]
\item $v_n = \sum_{d|n} \frac{d}{n} t_d b_{n/d}$, where the sum is over all integers $1<d\le n$ that divide $n$;
\item $\beta_n = -\sum_{s=1}^n \frac{s}{n} v_s \beta_{n-s}$;
\item for every positive integer $R$
\[
\sum_{s=0}^{R-1}\beta_sT_{n-s}=v_n+\frac{1}{n}\sum_{r=0}^{R-1}\beta_r\sum_{s=R-r}^{n-R} sv_sT_{n-r-s}\,.
\]

\end{enumerate}
\end{lem}

\noindent\emph{Proof:} Using formal term-by-term differentiation it is easy to show that \eqref{eq:1} and \eqref{eq:2} are equivalent to the formal identities
\begin{equation}\label{eq:3.1}
\log \left(\sum_{n=0}^\infty T_n x^n \right) = \sum_{m=1}^\infty v_m x^m
\end{equation}
and, respectively,
\begin{equation}\label{eq:3.2}
\log\left(\sum_{n=0}^\infty a_n x^n\right) = \sum_{s=1}^\infty \frac{b_s}{s}x^s\,.
\end{equation}
Taking the formal logarithm of \eqref{eq:0} and substituting \eqref{eq:3.1}, \eqref{eq:3.2}, we obtain
\[
\sum_{m=1}^\infty v_m x^m =\sum_{r,s=1}^\infty t_r \frac{b_s}{s} x^{rs} 
\]
from which 1) easily follows. Since \eqref{eq:2} is equivalent to the identity 
\[
\left(\sum_{s=0}^\infty \beta_s x^s\right)\left(\sum_{n=0}^\infty T_n x^n\right) = 1 
\]
of formal power series, then taking formal logarithms yields
\[
\log\left(\sum_{s=0}^\infty \beta_s x^s\right) = - \sum_{m=1}^\infty v_m x^m\,.
\]
Comparing with \eqref{eq:1} and \eqref{eq:3.1} proves 2). It follows from 2) that
\[
\sum_{u=0}^{R-1}\sum_{r=0}^u \beta_r \left( (n-u)v_{n-u}T_{u-r}+(u-r)v_{u-r}T_{n-u}\right) = v_n - \sum_{u=0}^{R-1} u \beta_u T_{n-u}
\]
and thus
\begin{align*}
\sum_{s=0}^{R-1} \beta_s T_{n-s} - v_n &= \frac{1}{n}\sum_{u=0}^{R-1}\sum_{r=0}^u \beta_r\left(\frac{n-r}{R-r}T_{n-r}- (n-u)v_{n-u}T_{u-r}-(u-r)v_{u-r}T_{n-u} \right)\\
&=\frac{1}{n}\sum_{r=0}^{R-1}\beta_r\left((n-r)T_{n-r}-\sum_{s=0}^{R-r-1}(n-r-s)T_sv_{n-r-s}+\sum_{s=0}^{R-r-1} s v_s T_{n-r-s} \right)\\
&= \frac{1}{n}\sum_{r=0}^{R-1}\beta_r\sum_{s=R-r}^{n-R} sv_sT_{n-r-s}
\end{align*}
where the last equality is obtained applying \eqref{eq:1} to $T_{n-r}$ for each $r\in \{0,1,\ldots,R-1\}$.

\begin{lem}\label{lem:3}
Let $(\{T_n\},\{t_n\},\{a_n\})$ be a Warlimont triple. Then 
\[
a_1\sum_{s=0}^{\lfloor n/2\rfloor} T_st_{n-s} \le T_n
\]
for all $n$.
\end{lem}

\noindent\emph{Proof:} Since
\[
\prod_{m=N+1}^\infty \left( \sum_{k=0}^\infty a_k x^{km}\right)^{t_m}\in 1+x^{N+1}{\mathbb R}[[x]]
\]
for every integer $N\ge 0$, then
\[
\sum_{n=0}^N T_n x^n  = \prod_{m=1}^N \left(\sum_{k=0}^\infty a_k x^{km}\right)^{t_m} + x^{N+1}{\mathbb R}[[x]]
\]
and thus
\begin{equation}\label{eq:7}
\sum_{n=0}^N T_n x^n  = \left(\sum_{s=0}^{\lfloor N/2\rfloor} T_s x^s \right)\prod_{m=\lfloor N/2\rfloor +1}^N \left(\sum_{k=0}^\infty a_k x^{km}\right)^{t_m} + x^{N+1}{\mathbb R}[[x]]\,.
\end{equation}
On the other hand by assumption $t_m$ is a non-negative integer for all $m$ and thus by the binomial theorem
\begin{equation}\label{eq:8}
\left(\sum_{k=0}^\infty a_k x^{km}\right)^{t_m} = 1 + a_1t_m x^m + x^{2m}{\mathbb R}[[x]]\,.
\end{equation}
Since the sequences $\{a_k\}$, $\{t_m\}$ and $\{T_n\}$ are non-negative, the Lemma follows by substituting \eqref{eq:8} into \eqref{eq:7} and comparing coefficients.

\begin{lem}\label{lem:4}
Let $(\{T_n\},\{t_n\},\{a_n\})$ be a Warlimont triple such that $\log(a_n)=O(n)$. Then for every non-negative integer $R$
\[
|v_n-a_1t_n|=
\begin{cases}
O(T_{n-R}) \quad \textrm{ if }\quad T_{n-1}=o(T_n)\\
\\
O(t_{n-R}) \quad \textrm{ if }\quad  t_{n-1}=o(t_n)\,.
\end{cases}
\]
\end{lem}

\noindent\emph{Proof:} Assume $T_{n-1}=o(T_n)$. Since $\log(a_n)=O(n)$, there exists $r>1$ such that $a_n\le r^n$ for all $n$. By induction of the definition of $\{b_n\}$, we obtain $|b_n|\le (3r)^n$ for all $n$. Moreover, since $T_{n-1}=o(T_n)$, there exists a constant $C$ such that $T_n\le C (3r)^{-2m} T_{n+m}$ for all $n,m$. Therefore, using Lemma 2 and Lemma 3
\[
|v_n-a_1t_n|\le \sum_{d/n} \frac{d}{n}t_d|b_{n/d}|\le \sum_{d/n} T_d (3r)^{n/d} \le CT_{n-R} \sum_{d/n}(3r)^{-n+2R+2d}=O(T_{n-R})\,. 
\]
The proof for the case $t_{n-1}=o(t_n)$ is similar and left to the reader.

\begin{thm}\label{thm:5}
Let $(\{T_n\},\{t_n\},\{a_n\})$ be a Warlimont triple such that $\log(a_n)=O(n)$ and let $R$ be a positive integer. Then the following are equivalent
\begin{enumerate}[1)]
\item $T_{n-1}=o(T_n)$ and
\[
\sum_{s=R}^{n-R} T_s T_{n-s} = O(T_{n-R})\,;
\]
\item $T_{n-1}=o(T_n)$ and
\[
a_1t_n=\sum_{s=0}^{R-1}\beta_s T_{n-s} + O(T_{n-R});
\]
\item $t_{n-1}=o(t_n)$ and
\[
T_n=a_1\sum_{s=0}^{R-1}T_st_{n-s}+O(t_{n-R});
\]
\item $t_{n-1}=o(t_n)$ and
\[
\sum_{s=R}^{n-R}t_st_{n-s}=O(t_{n-R})\,.
\]
\end{enumerate}
\end{thm}

\noindent\emph{Proof:} Assume 1) holds. Using Lemma \ref{lem:3}, Lemma \ref{lem:4} and $T_{n-1}=o(T_n)$, we obtain
\[
|v_n|\le |v_n-a_1t_n|+a_1t_n = O(T_n)\,.
\]
Therefore, exists an integer $N>R$ and a constant $C>0$ such that $|v_n|\le C T_n\le C T_{n+r}$ for all $n\ge N$ and for all $r\in \{0,\ldots,R-1\}$. Combining this observation with Lemma \ref{lem:2} yields
\begin{align*}
\left|v_n-\sum_{s=0}^{R-1}\beta_sT_{n-s}\right| & \le \sum_{r=0}^{R-1}\beta_r\sum_{s=R-r}^{n-R} \frac{s}{n}|v_s|T_{n-s}\\
& \le \sum_{r=0}^{R-1} \beta_r\left( \sum_{s=R-r}^{N-1} |v_s|T_{n-r-s}+C\sum_{s=N}^{n-R} T_sT_{n-r-s}\right)\\
&\le \sum_{r=0}^{R-1} \beta_r\left( \sum_{s=R-r}^{N-1} |v_s|T_{n-r-s}+C^2\sum_{s=R}^{n-R} T_sT_{n-s}\right)\\
& = O(T_{n-R})
\end{align*}
and thus
\[
\left|a_1t_n-\sum_{s=0}^{R-1}\beta_sT_{n-s}\right|\le |a_1t_n-v_n|+\left|v_n-\sum_{s=0}^{R-1}\beta_s T_{n-s}\right| = O(T_{n-R})\,.
\]
Hence, 1) implies 2). Assume 2) holds. Then
\begin{align*}
a_1\sum_{s=0}^{R-1}T_st_{n-s}&=\sum_{s=0}^{R-1}T_s\sum_{r=0}^{R-1-s}\beta_rT_{n-s-r}+O(T_{n-R})\\
&=\sum_{s=0}^{R-1}\sum_{u=s}^{R-1}T_s\beta_{u-s} T_{n-u}+O(T_{n-R})\\
&=\sum_{u=0}^{R-1}T_{n-u}\sum_{s=0}^uT_s\beta_{u-s}+O(T_{n-R}) \\
&=T_n+O(T_{n-R})
\end{align*}
where the last equality is obtained using the definition of the sequence $\{\beta_n\}$. In particular, setting $R=1$ we obtain $a_1t_n-T_n=O(T_{n-1})=o(T_n)$ which implies $a_1t_n\sim T_n$ and hence $o(t_{n-1})=O(T_{n-1})=o(T_n)=o(t_n)$. Therefore, 2) implies 3). Assume 3) holds. By Lemma \ref{lem:3}
\[
a_1^2\sum_{s=R}^{\lfloor n/2\rfloor} t_st_{n-s}\le a_1\sum_{s=R}^{\lfloor n/2\rfloor}T_st_{n-s}\le T_n-a_1\sum_{s=0}^{R-1}T_st_{n-s}=O(t_{n-R})\,.
\]
This proves 4) since 
\[
\sum_{s=R}^{n-R}t_st_{n-s}= 2 \sum_{s=R}^{\lfloor n/2\rfloor} t_st_{n-s} + O(t_{n-R})\,.
\]
Finally, assume 4) holds. Lemma \ref{lem:4} implies $|v_n-a_1t_n|=O(t_{n-R})=o(t_n)$ and thus $v_n\sim a_1t_n$. This implies that there exists an integer $N\ge R$ and constants $c,C>0$ such that
\begin{equation}\label{eq:4}
0<ct_n\le v_n\le Ct_n
\end{equation}
for all $n\ge N$. As a consequence,
\begin{equation}\label{eq:5}
v_{n-1}=O(t_{n-1})=o(t_n)=o(v_n)
\end{equation}
and 
\begin{equation}\label{eq:6}
\sum_{j=N}^{n-N}v_{n-j}v_j\le C^2\sum_{j=R}^{n-R} t_{n-j}t_j = O(t_{n-R})=O(v_{n-R})\,.
\end{equation}
For each $n\ge N$, let
\[
M_n=\max\left\{\frac{T_j}{v_j}\,|\, N\le j\le n \right\}\,.
\]
By Lemma \ref{lem:2}, we obtain
\begin{align*}
|T_n-v_n| &\le \sum_{j=1}^{n-1} |v_{n-j}|T_j\\ 
& \le\sum_{j=1}^{N-1}v_{n-j}T_j+M_{n-1}\left(\sum_{j=N}^{n-N} v_{n-j}v_j +\sum_{j=n-N+1}^{n-1}|v_{n-j}|v_j\right)\\
&= o(v_n)(1+M_{n-1})
\end{align*}
where \eqref{eq:5} and \eqref{eq:6} where used to obtain the last equality. Hence there exists $N_1\ge N$ such that for all $n\ge N_1$
\[
\frac{T_n}{v_n} \le \frac{3+M_{n-1}}{2}
\]
and thus
\[
M_n=\max\left\{M_{n-1},\frac{T_n}{v_n}\right\}\le \max\{M_{n-1},3\}\,.
\]
This shows that the sequence $\{M_n\}$ is bounded i.e.\ there exists a constant $K>0$ such that $T_n\le K v_n$ for all $n\ge N$. Therefore, using \eqref{eq:4} and Lemma \ref{lem:3}, we obtain
\[
T_{n-1}=O(v_{n-1})=O(t_{n-1})=o(t_n)=o(T_n)\,.
\]
Moreover \eqref{eq:5} yields 
\[
\sum_{s=R}^{n-R} T_{n-s}T_s \le 2\sum_{s=R}^N T_{s}T_{n-s} + K^2 \sum_{s=N}^{n-N} v_{n-s}v_s = O(T_{n-R})+O(v_{n-R}) = O(T_{n-R})\,. 
\]
This concludes the proof that 4) implies 1) and the Theorem is proved.

\begin{rem}\label{rem:6}
Let $(\{T_n\},\{t_n\},\{a_n\})$ be a Warlimont triple that satisfies the equivalent conditions of Theorem \ref{thm:5} for some $R>2$. Then
\[
\sum_{s=R-1}^{n-R}T_sT_{n-s}=\sum_{s=R}^{n-R}T_sT_{n-s}+2T_{R-1}T_{n-R+1}=O(T_{n-R})+O(T_{n-R+1})=O(T_{n-R+1})
\]
and thus $(\{T_n\},\{t_n\},\{a_n\})$ satisfies the equivalent conditions of Theorem \ref{thm:5} for any positive integer less or equal than $R$. In particular, $t_{n-1}=o(t_n)$ and $T_n\sim a_1t_n$.
\end{rem}

\section{Warlimont Functions and Wright Semigroups}\label{sec:3}

\begin{dfn}\label{def:6}
An {\it additive arithmetical semigroup} is a pair $(G,+,\partial)$ consisting of an abelian semigroup $(G,+)$ with identity and a semigroup homomorphism $\partial:(G,+)\to (\mathbb Z_{\ge 0},+)$ such that 
\begin{enumerate}[i)]
\item the cardinality $G_n$ of the preimage $\partial^{-1}(n)$ is finite for all $n$;
\item $G$ is freely generated by $G^+\subseteq G$.
\end{enumerate}
We denote by $G_n^+$ the cardinality of the set $\partial^{-1}(n)\cap G^+$.
\end{dfn}

\begin{rem}\label{rem:8}
Let $(G,+,\partial)$ be an additive arithmetical semigroup. As pointed out in \cites{K,warlimont}, $(\{G_n\},\{G_n^+\},\{1\})$ is a Warlimont triple.
\end{rem}

\begin{dfn}\label{def:7}
A {\it Wright semigroup} is an additive arithmetical semigroup $(G,+,\partial)$ such that 
\begin{equation}\label{eq:11.01}
\log(G_n)=\alpha n^{a+1} +\beta n\log(n) + \gamma n + O(n^b)
\end{equation}
for some real numbers $\alpha,\beta,\gamma,a,b$ such that $\alpha>0$ and $0<b<a$.
\end{dfn}

\begin{dfn}\label{def:10}
Let $R$ be a positive integer. We say that an additive arithmetical semigroup $(G,+,\partial)$ satisfies {\it axiom $\mathcal W_R$} if $G_{n-1}=o(G_n)$ and
\[
\sum_{s=R}^{n-R} G_s G_{n-s} = O(G_{n-R})\,.
\]
\end{dfn}

\begin{rem}\label{rem:11}
Let $(G,+,\partial)$ be an additive arithmetical semigroup that satisfies axiom $\mathcal W_R$ for some positive integer $R$. Combining Remark \ref{rem:6} and \ref{rem:8} we conclude that $(G,+,\partial)$ satisfies axiom $\mathcal W_{R'}$ for any positive integer $R'\le R$. In particular, $G_n\sim G_n^+$ and $G_{n-1}^+=o(G_n^+)$ i.e.\ the additive arithmetical semigroup $(G,+,\partial)$ satisfies both axiom $\mathcal G_1$ and axiom $\mathcal G_2$ as defined in \cite{K}. Notice that the combination of Axiom $\mathcal G_1$ and Axiom $\mathcal G_2$ is slightly weaker than axiom $\mathcal W_1$ since $\sum_{s=1}^{n-1}G_sG_{n-s}=o(G_n)$ does not necessarily imply $\sum_{s=1}^{n-1}G_sG_{n-s}=O(G_{n-1})$.
\end{rem}

\begin{prop}\label{prop:11}
Every Wright semigroup satisfies axiom $\mathcal W_R$ for every positive integer $R$.
\end{prop}

\noindent\emph{Proof:} This is a straightforward consequence of the definitions and Theorem 7 of \cite{wright70}.

\begin{dfn}\label{def:12}
Let $(G,+,\partial)$ be an additive arithmetical semigroup. A function $F:G\to \mathbb R$ is {\it multiplicative} if
$F(g_1+g_2)=F(g_1)F(g_2)$ for all $g_1,g_2\in G$ coprime. We say that $F$ is {\it prime-independent} if there exists a sequence $\{F^+_n\}$ such that $F^+_n=F(np)$ for every $p\in G^+$ and every positive integer $n$. For every function $F:G\to \mathbb R$, we denote by $\{F_n\}$ the sequence defined by setting
\[
F_n = \sum_{\partial(g)=n} F(g) 
\]
for each non-negative integer $n$. A {\it Warlimont function} is a non-negative multiplicative prime-independent function such that $\log(F_n^+)=O(n)$ and $F_1^+>0$. The {\it normalization} of a Warlimont function $F$ is the (not necessarily multiplicative) function $\widetilde F:G\to \mathbb R$ such that $\widetilde F(g)=F(g)/F_1^+$ for all $g\in G$.
\end{dfn}

\begin{ex}\label{ex:11}
Let $(G,+,\partial)$ be an additive arithmetical semigroup and let $F:G\to \mathbb R$ be such that $F(g)=1$ for all $g\in G$. Then $F$ is a Warlimont function and $F_n=\widetilde F_n=G_n$ for all $n$. 
\end{ex}

\begin{ex}\label{ex:12}
Let $(G,+,\partial)$ be an additive arithmetical semigroup and, for each $k\ge 2$, consider the {\it generalized divisor function} $d_k:G\to \mathbb R$ that to each $g\in G$ assigns the number $d_k(g)$ of $k$-tuples $(g_1,\ldots,g_k)\in G^k$ such that $g=g_1+\ldots+g_k$. Then $d_k$ is multiplicative, prime-independent and $(d_k)^+_n=\binom{n+k-1}{k-1}$ for each integer $n\ge 1$. Therefore, $d_k$ is Warlimont.
\end{ex}

\begin{ex}\label{ex:13}
Let $(G,+,\partial)$ be an additive arithmetical semigroup and consider the {\it unitary divisor function} $d_*:G\to \mathbb R$ that to each $g\in G$ assigns the number $d_*(g)$ of coprime pairs $(g_1,g_2)$ such that $g=g_1+g_2$. Then $d_*$ is multiplicative, prime-independent and $(d_*)^+_n=2$ for each integer $n\ge 1$. Therefore $d_*$ is Warlimont.
\end{ex}

\begin{ex}\label{ex:14}
Let $(G,+,\partial)$ be an additive arithmetical semigroup and consider the {\it prime divisor function} $B:G\to \mathbb R$ such that
$B(k_1p_1+k_2p_2\ldots+k_rp_r)=k_1k_2\cdots k_r$ for any $p_1,\ldots,p_r\in G$ primes and $k_1,\ldots,k_r$ positive integers. Then $B$ is multiplicative, prime-independent and $B^+_n=n$ for each integer $n\ge 1$. Therefore, $B$ is Warlimont.
\end{ex}

\begin{rem}\label{rem:14}
Let $F$ be a Warlimont function on an additive arithmetical semigroup $(G,+,\partial)$. Then $F^m:G\to \mathbb R$ such that $F^m(g)=(F(g))^m$ is again a Warlimont function for every integer $m\ge 1$ since
\[
\log((F^m)^+_n) = m\log(F^+_n) = O(n)\,. 
\]
Moreover, $\widetilde{F^m}=(\widetilde F)^m$.
\end{rem}

\begin{rem}\label{rem:15}
Let $F$ be a Warlimont function on an additive arithmetical semigroup $(G,+,\partial)$. Then, as observed in \cite{warlimont}, $(\{F_n\},\{G_n^+\},\{F^+_n\})$ is a Warlimont triple.
\end{rem}

\begin{thm}\label{thm:16}
Let $(G,+,\partial)$ be an additive arithmetical semigroup that satisfies axiom $\mathcal W_R$ and let $F$ be a Warlimont function on $G$. Then for every positive integer $M$ there exist constants $\xi_1,\ldots \xi_{R-1}$ such that
\begin{equation}\label{eq:11.1}
\sum_{\partial(g)=n}(\widetilde{F}(g)-1)^M= \sum_{s=1}^{R-1} \xi_s G_{n-s} + O (G_{n-R})\,.
\end{equation}

\end{thm}

\noindent\emph{Proof:} $(\{G_n\},\{G_n^+\},\{1\})$ and $(\{F_n\},\{G_n^+\},\{F_n^+\})$ are both Warlimont triple by Remark \ref{rem:15} and Example \ref{ex:11}. Since $\{G_n\}$ satisfies axiom $\mathcal W_R$, it follows from Theorem \ref{thm:5} applied to the Warlimont triple $(\{G_n\},\{G_n^+\},\{1\})$ that $G^+_{n-1}=o(G_n^+)$,
\begin{equation}\label{eq:12}
\sum_{s=R}^{n-R} G_s^+G_{n-s}^+= O (G_{n-R}^+)\,.
\end{equation}
Moreover, if $\{\beta_n\}$ is the sequence defined recursively by setting $\beta_0=1$ and
\begin{equation}\label{eq:14.1}
\beta_n=-\sum_{s=0}^{n-1}\beta_sG_{n-s}
\end{equation}
for every positive integer $n$, then 
\begin{equation}\label{eq:13}
G_{n-s}^+=\sum_{r=0}^{R-1}\beta_r G_{n-s-r}+O(G_{n-s-R})\,.
\end{equation}
for all $s\ge 0$. In particular, we can apply Theorem \ref{thm:5} to the Warlimont triple $(\{F_n\},\{G_n^+\},\{F_n^+\})$  and obtain
\begin{equation}\label{eq:14}
F_n=F_1^+\sum_{s=0}^{R-1} F_s G_{n-s}^+ + O(G_{n-R}^+)\,.
\end{equation}
Since by definition $G_n^+\le G_n$ for all $n$ and $G_{n-s-R}=o(G_{n-R})$ for all $s>0$, then substituting \eqref{eq:13} into \eqref{eq:14} yields
\begin{equation}\label{eq:15.1}
\widetilde F_n = \sum_{s=0}^{R-1}\left(\sum_{r=0}^s \beta_r F_{s-r}\right)G_{n-s}+O(G_{n-R})
\end{equation} 
Using the binomial theorem and Remark \ref{rem:14} we obtain
\begin{align}
\sum_{\partial(g)=n}(\widetilde F(g)-1)^M & = (-1)^M \sum_{\partial(g)=n}\sum_{m=0}^M (-1)^m \binom{M}{m}\widetilde{F^m}(g) \\ &= (-1)^M\sum_{m=0}^M (-1)^m \binom{M}{m} (\widetilde{F^m})_n\,. \label{eq:15}
\end{align}
Applying \eqref{eq:15.1} to the Warlimont function $F^m$ and substituting into \eqref{eq:15} (after an obvious rearrangement) yields
\[
\sum_{\partial(g)=n}(\widetilde F(g)-1)^M= \sum_{s=0}^{R-1} \xi_s G_{n-s} + O (G_{n-R})\,,
\]
with
\begin{align}
\xi_s &= (-1)^M \sum_{m=0}^M(-1)^m\binom{M}{m}\sum_{r=0}^s\beta_r (F^m)_{s-r}\\ &= (-1)^M \sum_{m=1}^M(-1)^m\binom{M}{m}\sum_{r=0}^s\beta_r (F^m)_{s-r} \label{eq:17.1}
\end{align}
for all $s\in\{0,\ldots,R-1\}$ where the second equality follows from \eqref{eq:14.1} and Example \ref{ex:11}.
This implies \eqref{eq:11.1} since (by combining Remark \ref{rem:14} and Remark \ref{rem:15}) $(F^m)_0=1$ for all $m$ and thus
\[
\xi_0= (-1)^M\sum_{m=0}^M(-1)^m\binom{M}{m} \beta_0 (F^m)_0 =0\,.
\]

\begin{dfn}
Let $F$ be a Warlimont function on an additive arithmetical semigroup $(G,+,\partial)$ and let $M$ be a positive integer. We define the {\it normalized $M$-th moments of $F$} to be the functions $\mu_{F,M}:\mathbb Z_{\ge 0}\to \mathbb R$ defined by
\[
\mu_{F,M}(n)=\frac{1}{G_n}\sum_{\partial(g)=n} (\widetilde{F}(g)-1)^M
\]
for all $n\ge 0$. 
\end{dfn}

\begin{rem}\label{rem:22}
Let $F$ be a Warlimont function on an additive arithmetical semigroup $(G,+,\partial)$. The average value of $F$ on $\partial^{-1}(n)$ is given by
\[
\frac{F_n}{G_n}=F_1^+(1+\mu_{F,1}(n))\,.
\] 
The higher normalized moments can be thought of as capturing the deviation of $F$ from $F_1^+$. For instance, if $\mu_{F,1}(n)=o(1)$, then 
\[
\frac{1}{G_n}\sum_{\partial(g)=n}(F(g)-F_1^+)^2 = (F_1^+)^2\mu_{F,2}(n)
\]
can be thought of as an asymptotic measure of the variance of $F$ on $\partial^{-1}(n)$.
\end{rem}

\begin{cor}\label{cor:23}
Let $F$ be a Warlimont function on an additive arithmetical semigroup $(G,+,\partial)$ that satisfies axiom $W_1$. Then
\[
\lim_{n\to \infty} \frac{F_n}{G_n} = F_1^+
\]
and
\[
\lim_{n\to \infty} \frac{1}{G_n} \sum_{\partial(g)=n} (F(g)-F_1^+)^2 =0\,.
\]
\end{cor}

\noindent\emph{Proof:} Combining Remark \ref{rem:22} and Theorem \ref{thm:16} (with $R=M=1$), we obtain
\[
\frac{F_n}{G_n} = F_1^+(1+\mu_{F,1}(n))= F_1^+ + o(1)\,.
\]
Similarly,
\[
\frac{1}{G_n} \sum_{\partial(g)=n} (F(g)-F_1^+)^2 = (F_1^+)^2 \left(\xi_1 \frac{G_{n-1}}{G_n} +O\left(\frac{G_{n-1}}{G_n}\right)\right) = o(1)\,.
\]

\begin{rem}
A slightly stronger (see Remark \ref{rem:11}) version of Corollary \ref{cor:23} is proved in \cite{K} for particular choices of $F$. A sharper result is given in \cite{warlimont} where it is shown that the assumption $G_{n-1}=O(G_n)$ (which is part of axiom $\mathcal W_1$) is unnecessary.
\end{rem}

\begin{thm}\label{thm:25}
Let $F$ be a Warlimont function on a Wright semigroup $(G,+,\partial)$ with $\alpha,a,b$ as in Definition \ref{def:7} and let $q=e^{\alpha(a+1)}$.  
\begin{enumerate}[1)]
\item For every positive integer $M$ there exists a sequence $\{\lambda_s\}$ of functions $\lambda_s:\mathbb Z_{\ge 0}\to \mathbb R$ such that $\log(\lambda_s(n))=O(n^{a-1}+n^b)$ and
\begin{equation}\label{eq:18.1}
\mu_{F,M}(n)=\sum_{s=1}^{R-1} \lambda_s(n) q^{-sn^a} + O\left(\lambda_R(n)q^{-Rn^a}\right)
\end{equation}
for every integer $R>0$.
\item Assume further that there exists constants $0\le d_2\le d_1$ and sequence $\{\psi_s\}$ of polynomials such that $\deg(\psi_s)\le d_1s-d_2$ for all $s\ge 1$ and 
\begin{equation}\label{eq:19}
\frac{G_{n-1}}{G_n} = \sum_{s=1}^{R-1} \psi_s(n) q^{-ns} +O(n^{d_1R-d_2}q^{-Rn})
\end{equation}
for every integer $R>0$. Then there exists a sequence $\{\tau_s\}$ of polynomials such that $\deg(\tau_s)\le d_1s-d_2$ and
\begin{equation}\label{eq:20}
\mu_{F,M}(n)=\sum_{s=1}^{R-1} \tau_s(n) q^{-sn} + O\left(n^{d_1R-d_2}q^{-Rn}\right)
\end{equation}
for every integer $R>0$.
\end{enumerate}

\end{thm}

\noindent\emph{Proof:} Let $\xi_s$ be defined by \eqref{eq:17.1} for all $s\ge 1$. By Proposition \ref{prop:11}, and Theorem \ref{thm:16},
\begin{equation}\label{eq:22}
\mu_{F,M}(n)=\sum_{s=1}^{R-1} \xi_s \frac{G_{n-s}}{G_n} + O \left(\frac{G_{n-R}}{G_n}\right)
\end{equation}
for every integer $R>0$. Since
\[
\log\left(\frac{G_{n-s}}{G_n} \right)=\alpha((n-s)^{a+1}-n^{a+1}) +O (n^b) = -\alpha(a+1)sn^a + O(n^{a-1}+n^b)\,,
\]
then in order to prove 1) it suffices to choose $\lambda_s$ such that
\[
\lambda_s(n)= \xi_s q^{sn^a}\frac{G_{n-s}}{G_n}
\]
for all $n\ge s \ge 1$. Using \eqref{eq:19} repeatedly and induction on $t$, we obtain
\begin{equation}\label{eq:23}
\frac{G_{n-t}}{G_n}=\frac{G_{n-1}}{G_n}\cdots\frac{G_{n-t}}{G_{n-t+1}}=\sum_{s=t}^{R-1} \nu_{s,t}(n) q^{-sn} + O(n^{DR}q^{-Rn})
\end{equation}
where
\begin{equation}\label{eq:25}
\nu_{s,t}(n)=\sum_{i_1+\cdots+i_t=s}  \psi_{i_1}(n)\psi_{i_2}(n-1)\cdots \psi_{i_t}(n-t+1) q^{i_2+2i_3+\cdots + (t-1) i_t}\,.
\end{equation}
is a polynomial in $n$ of degree at most $d_1s-d_2t$ for all $1\le t \le s$. Substituting \eqref{eq:23} in into \eqref{eq:22} yields
\begin{align*}
\mu_{F,M}(n) & = \sum_{t=1}^{R-1}\xi_t\sum_{s=t}^{R-1}\nu_{s,t}(n)q^{-ns}+O(n^{DR}q^{-nR})\\ 
& = \sum_{s=1}^{R-1} \left(\sum_{t=1}^s \xi_t\nu_{s,t}(n) \right) q^{-sn} + O(n^{DR}q^{-nR})\,,
\end{align*}
which proves 2) upon setting
\begin{equation}\label{eq:28.1}
\tau_s(n)=\sum_{t=1}^s \xi_t\nu_{s,t}(n)
\end{equation}
for all $s,n$.

\begin{rem}\label{rem:26}
Comparison of \eqref{eq:19} and \eqref{eq:11.01} shows that the assumptions of 2) in Theorem \ref{thm:25} require in particular that \eqref{eq:11.01} holds with $a=1$.
\end{rem}

\section{Examples}\label{sec:4}

\subsection{Graphs}\label{sec:4.1} Let $(G,+)$ be the semigroup of (simple, unlabeled) graphs with semigroup operation $+$ given by disjoint union. If $\partial$ is the map that to each graph $g$ assigns the cardinality of its set of vertices, then $(G,+,\partial)$ is an additive arithmetical semigroup and $g\in G^+$ if and only if the graph $g$ is connected. As proved in \cite{wright69}, there exists a sequence $\{\varphi_s\}$ of polynomials such that $\varphi_s$ has degree $2s$ for every $s$ and
\begin{equation}\label{eq:26}
G_n = \frac{2^{\binom{n}{2}}}{n!}\left(\sum_{s=0}^{R-1}\varphi_s(n) 2^{-sn} + O(n^{2R}2^{-Rn}) \right)
\end{equation}
for every positive integer $R$. The polynomials $\varphi_s$ can be calculated explicitly, the first few being
\begin{eqnarray*}
\varphi_0(n)&=& 1\,;\\
\varphi_1(n)&=&2n^2-2n\,;\\
\varphi_2(n)&=& 8 n^4 - \frac{128}{3} n^3 + 72 n^2 - \frac{112}{3} n\,;\\
\varphi_3(n)&=&\frac{256}{3} n^6 -\frac{3712}{3} n^5 + \frac{20672}{3} n^4 -\frac{54272}{3} n^3 + 21952 n^2 - 9600 n\,.
\end{eqnarray*}
In particular, 
\[
\log(G_n)=\log(\sqrt{2})n^2-n\log(n)+ (1-\log(\sqrt{2}))n+O(n^b)
\]
for any $b>0$ and thus $(G,+,\partial)$ is a Wright semigroup. Moreover, using \eqref{eq:26} and expanding the denominator as a geometric series we obtain
\begin{align*}
\frac{G_{n-1}}{G_n} &=2n2^{-n}\left(\sum_{s=0}^{R-1}2^s\varphi_s(n-1)2^{-sn}\right)\sum_{r=0}^{R-1}\left(-\sum_{s=1}^{R-1}\varphi_s(n)2^{-sn}\right)^r + O(n^{2R+1}2^{-(R+1)s})\\
&=\sum_{s=1}^{R-1}\psi_s(n)2^{-sn}+O(n^{2R-1}2^{-Rn})\,,
\end{align*}
where the $\psi_s$ are polynomials of degree $\deg(\psi_s)=2s-1$ which can be explicitly calculated in terms of the polynomials $\varphi_s$ in \eqref{eq:26}. For instance
\begin{eqnarray*}
\psi_1(n)&=&2n\,;\\
\psi_2(n)&=&4n^3-20n^2+16n\,;\\
\psi_3(n)&=&40 n^5 - 464 n^4  + 1768 n^3 - 2624 n^2+1280 n  \,;\\
\psi_4(n) &=&\frac{3248}{3} n^7- 24176 n^6  + \frac{630608}{3}n^5 - 908496 n^4 + \frac{6137792}{3} n^3 - 2250240 n^2 + 925696 n     \, .
\end{eqnarray*}
Substitution into \eqref{eq:25} yields $\nu_{s,1}(n)=\psi_s(n)$ for all $s$ and
\begin{eqnarray*}
\nu_{2,2}(n)&=& 8 n^2 - 8n \,;\\
\nu_{3,2}(n)&=& 48n^4-352n^3+688n^2-384n\,;\\
\nu_{3,3}(n)&=& 64n^3-192n^2+128n\,;\\
\nu_{4,2}(n)&=&  864 n^6 - 13472 n^5 + 77216 n^4 - 203488 n^3  + 245376 n^2  -106496 n  \,;\\
\nu_{4,3}(n)&=& 896n^5-9728n^4+35200n^3-50944n^2+24576n\,;\\
\nu_{4,4}(n)&=& 1024n^4-6144n^3+11264n^2-6144n\,.
\end{eqnarray*}
Inspection of graphs with up to four vertices shows that $G_1=1$, $G_2=2$, $G_3=4$ and $G_4=11$. Substitution into \eqref{eq:14.1} yields $\beta_1=\beta_2=\beta_3=-1$ and $\beta_4=-4$. 

\begin{ex}
Consider the Warlimont function $d_2$ from Example \ref{ex:12}. When specialized to the semigroup of graphs, $d_2$ counts the number of ways of writing a given graph as the disjoint union of two graph. The order is taken into account, so that if $g_1$ is not isomorphic to $g_2$, then $g=g_1+g_2$ and $g=g_2+g_1$ count as two distinct decompositions. Moreover, decompositions in which one of the components is the empty graph are allowed. Combining Remark \ref{rem:22} and Theorem \ref{thm:25} we obtain \eqref{eq:00000}. In particular, setting $M=1$ yields a full asymptotic expansion for the average of $d_2$ of the form
\[
\frac{1}{G_n}\sum_{\partial(g)=n} d_2(g) = 2 + 2\sum_{s=1}^{R-1}\tau_s(n) 2^{-sn}+O\left(n^{2R-1}2^{-Rn}\right),
\]
valid for every positive integer $R$ where the $\tau_s(n)$ are polynomials of degree $2s-1$. For instance, direct inspection of graphs with up to four vertices yields $(d_2)_1=2$, $(d_2)_2=5$, $(d_2)_3=12$ and $(d_2)_4=34$. Substituting into \eqref{eq:17.1} and then into \eqref{eq:28.1} we obtain
\begin{align*}
\tau_1(n) & = 2n \,;\\
\tau_2(n) & = 4 n^3-4n^2\,;\\
\tau_3(n) & = 40 n^5  - 368 n^4  + 1320 n^3 - 2016 n^2 + 1024 n \,;\\
\tau_4(n) &= \frac{3248}{3}n^7 - 22448 n^6 + \frac{560528}{3}n^5 - 781712 n^4 + \frac{5136512}{3}n^3 - 1839360 n^2  +743424 n    \,.
\end{align*}
\end{ex}

\subsection{Graphs with an even number of edges} Let $(G,+)$ be the semigroup of (simple, unlabeled) graphs with an even number of edges with semigroup operation $+$ given by disjoint union. If $\partial$ is the map that to each graph $g$ assigns the cardinality of its set of vertices, then $(G,+,\partial)$ is an additive arithmetical semigroup. $G^+$ consists of graphs $g$ with an even number of edges that cannot be written as the disjoint union of two nonempty graphs with an even number of edges. While $G$ is a subsemigroup of the semigroup of all unlabeled graphs, not all graphs in $G^+$ are connected. For instance, while $2K_1$ is not connected it is nevertheless prime in the semigroup of graphs with even edges. As pointed out in \cite{aldi}, there exists a sequence $\{\varphi_s\}$ of polynomials such that $\varphi_s$ has degree $2s$ for every $s$ and
\begin{equation*}
G_n = \frac{2^{\binom{n}{2}}}{2n!}\left(\sum_{s=0}^{R-1}\varphi_s(n) 2^{-sn} + O(n^{2R}2^{-Rn}) \right)
\end{equation*}
for every positive integer $R$, where the polynomials $\varphi_s(n)$ coincide with those of Section \ref{sec:4.1}.
In particular, $(G,+,\partial)$ is a Wright semigroup and
\[
\frac{G_{n-1}}{G_n} = \sum_{s=1}^{R-1}\psi_s(n)2^{-sn}+O(n^{2R-1}2^{-Rn})\,
\]
where the polynomials $\psi_s(n)$ coincide with those calculated in Section \ref{sec:4.1}. Inspection of graphs with up to four vertices shows that $G_1=G_2=1$, $G_3=2$ and $G_4=6$. Substitution into \eqref{eq:14.1} yields $\beta_1=-1$, $\beta_2=0$, $\beta_3=-1$ and $\beta_4=-3$. 

\begin{ex}
Consider the Warlimont function $d_*$ from Example \ref{ex:13}. Combining Remark \ref{rem:22} and Theorem \ref{thm:25} we obtain a full asymptotic expansion for the second moment of of $d_*$ about $2$ is
\[
\frac{1}{G_n}\sum_{\partial(g)=n} (d_*(g)-2)^2 = 4\sum_{s=1}^{R-1}\tau_s(n) 2^{-sn}+O\left(n^{2R-1}2^{-Rn}\right)
\]
for every positive integer $R$, where the $\tau_s(n)$ are polynomials of degree $2s-1$. To calculate these explicitly for small values of $s$, we first observe that (by direct calculation $(d_*)_1=2$, $(d_*)_2=2$, $(d_*)_3=4$, $(d_*)_4=14$ as well as $(d_*^2)_1=4$, $(d_*^2)_2=4$, $(d_*^2)_3=8$, $(d_*^2)_4=36$. Substitution into \eqref{eq:17.1} (upon setting $M=2$) and then into \eqref{eq:28.1} yields
\begin{align*}
\tau_1(n) & =2n  \,;\\
\tau_2(n) & = 4 n^3 - 20 n^2 + 16 n\,;\\
\tau_3(n) & = 40n^5 - 464 n^4 + 1832 n^3 - 2816 n^2 +1408n  \,;\\
\tau_4(n) & = \frac{3248}{3}n^7- 24176 n^6+ \frac{633296}{3}n^5 - 906960 n^4 + \frac{6040640}{3}n^3- 2177280 n^2+882688 n       \,.
\end{align*}
\end{ex}

\subsection{Polynomials over a finite field} Consider the field $\mathbb F_q$ with $q$ elements and let $G$ be the set of non-zero polynomials in $\mathbb F_q[x_1,\ldots,x_k]$ modulo the equivalence relation $f\sim g$ if and only if $f=\lambda g$ for some $\lambda\in \mathbb F_q$. $G$ has a natural structure of additive semigroup with semigroup operation $+$ given by multiplication of polynomials. If $\partial$ is the semigroup homomorphism that to each polynomial $f\in G$ assigns its total degree, then $(G,+,\partial)$ is an additive arithmetical semigroup and $G^+$ is the set of equivalent classes of irreducible polynomials in $\mathbb F_q[x_1,\ldots,x_k]$. Since
\begin{equation}\label{eq:31}
G_n=\frac{q^{\binom{n+k}{k}}-q^{\binom{n-1+k}{k}}}{q-1}\,
\end{equation}
for every $n$, then
\[
\log(G_n)=\log(q)\frac{n^k}{k!} + O(n^{k-1})
\]
for every $k\ge 2$. On the other hand if $k=1$, then $\log(G_n)=\log(q)n$ for every $n$. Hence $(G,\cdot,\partial)$ is a Wright semigroup if and only if $k\ge 2$. If $k=2$ then for every positive integer $R$
\[
\frac{G_{n-1}}{G_n}=q^{-n-1} \frac{1-q^{-n}}{1-q^{-n-1}} = \sum_{s=1}^{R-1} \psi_{s}(n) q^{-sn} + O(q^{-Rn})
\]
where, $\psi_1(n)=q^{-1}$ and $\psi_s(n)=q^{-s}(1-q)$ for all $s\ge 2$. By Theorem \ref{thm:25}, each $\mu_{F,M}$ admits an asymptotic expansion as a power series in $q^{-n}$ with {\it constant} coefficients. For instance, substitution into \eqref{eq:25} yields
\begin{eqnarray*}
\nu_{2,1}(n)&=& q^{-2} - q^{-1}\,; \\
\nu_{2,2}(n)&=& q^{-1} \,;\\
\nu_{3,1}(n)&=& q^{-3}-q^{-2}\,;\\
\nu_{3,2}(n)&=& q^{-2}-1\,;\\ 
\nu_{3,3}(n)&=& 1\,.
\end{eqnarray*}
\begin{ex}
Let us further specialize to the case where $G$ is the semigroup of non-zero polynomial in two variables over the field with two elements. By Theorem \ref{thm:25}, there exist constants $\tau_s$ such that for every positive integer $R$ the average of the Warlimont function $B$ (as  defined in Example \ref{ex:14}) on polynomials of degree $n$ is
\begin{equation}\label{eq:32}
\frac{B_n}{G_n}=1+\sum_{s=1}^{R-1}\tau_s 2^{-sn} + O(2^{-Rn})\,.
\end{equation}
Since $B_1=6$, $B_2=62$ and $B_3=1002$, substituting \eqref{eq:31} into \eqref{eq:14.1} and then into \eqref{eq:17.1} shows that in particular $\tau_1=0$, $\tau_2=3$ and $\tau_3=\frac{3}{2}$.
\end{ex}

\begin{ex}
If $k>2$, then by Remark \ref{rem:26} the second part of Theorem \ref{thm:25}. Nevertheless, the asymptotic behavior of Warlimont functions can be described using \eqref{eq:18.1} as follows. Consider for instance the Warlimont function $B$ of Example \ref{ex:14} on the semigroup of polynomials in $3$ variables with coefficients in $\mathbb F_q$. Since
\[
G_1=-\beta_1= (B^m)_1
\]
for all $m$, substitution in \eqref{eq:17.1} yields $\xi_1=0$ and thus
\[
\frac{1}{G_n}\sum_{\partial(g)=n}(B(g)-1)^M=O\left(\frac{G_{n-2}}{G_n}\right)=O\left(q^{-n^2-2n} \right)
\]
for all $M$.
\end{ex}

\end{document}